\documentclass[12pt]{amsart}
\usepackage{graphicx}
\usepackage{amssymb}
\newtheorem{theorem}{Theorem}[section]

\newtheorem{proposition}[theorem]{Proposition}
\newtheorem{corollary}[theorem]{Corollary}
\newtheorem{definition}{Definition}

\newtheorem{remark}{Remark}

\def\mod{\operatorname{mod}}
\def\Fol{\operatorname{Fol}}
\def\mes{\operatorname{mes}}

\def\nit{\mathbb Z}

\author{A.Ya.Belov, G.V.Kondakov, I.Mitrofanov}

\title{Inverse problems of symbolic dynamics}
\date{}
\begin{document}

\footnote{Shanghai University, Moscow Institute of open Education,
Moscow Institute of Physicks and Technology, Moscow State
 University}

\begin{abstract}
This paper reviews some results regarding symbolic dynamics,
correspondence between languages of dynamical systems and
combinatorics. Sturmian sequences provide a pattern for
investigation of one-dimensional systems, in particular interval
exchange transformation. Rauzy graphs language can express 
many important combinatorial and some dynamical properties. In this case combinatorial
properties are considered as being generated by substitutional system, and dynamical properties are considered as
criteria of superword being generated by interval exchange
transformation. As a
consequence, one can get a morphic word appearing in interval
exchange transformation such that frequencies of letters are
algebraic numbers of an arbitrary degree.

Concerning multydimensional systems, our main result is the following.
Let $P(n)$ be a polynomial, having an irrational coefficient of
the highest degree. A word $w$ $(w=(w_n), n\in \nit)$ consists of
a sequence of first binary numbers of $\{P(n)\}$ i.e.
$w_n=[2\{P(n)\}]$. Denote  the number of different
subwords of $w$ of length $k$ by $T(k)$ .

\medskip
{\bf Theorem.} {\it There exists a polynomial $Q(k)$, depending
only on the power of the polynomial $P$, such that $T(k)=Q(k)$ for
sufficiently great $k$. }
\end{abstract}

\maketitle

\section{Introduction}

Methods of symbolic dynamics are rather useful in the study of
combinatorial properties of words, investigation of problems of
number theory and theory of dynamical systems. Let $M$ be a
compact metric space, $U\subset M$ be its open subspace, $f:M\to
M$ be a homeomorphism of the compact into itself, and $x\in M$ be
an initial point. It determines a sequence of points
$$x,f(x),\ldots,f^{(n)}(x),\ldots$$
With the sequence of iterations, one can associate an infinite
binary word

$$
w_n=\left\{
\begin{array}{rcl}
   a,\ f^{(n)}(x_0)\in U\\
   b,\ f^{(n)}(x_0)\not\in U\\
\end{array}\right.
$$
which is called the {\it evolution} of point $x_0$. If $f$ is
invertible then $n\in \mathbb{Z}$, otherwise $n\in \mathbb{N}$.
Symbolic dynamics investigates the interrelation between the
properties of the dynamical system $(M,f)$ and the combinatorial
properties of the word $W_n$. For words over alphabets which
comprise more symbols, several characteristic sets should be
considered: $U_1,\ldots,U_n$. Technical notions regarding
combinatorics of words see in section \ref{SbScCmplUnirec}.

 {\it Direct problem} of
symbolic dynamics consists in the description of properties of the word $W$,
based on the information about dynamic system. {\it Reverse
problem} consists in the description of $(M,f,U,x)$, based on the information about
$W$.

We shall point out some facts, which are known from folklore.

Minimality of the dynamical system corresponds with the uniform recurrence
property (see section \ref{SbScCmplUnirec}).

Uniqueness of invariant measure corresponds with the following property.
Let $u$ be a subword of uniformly recurrent (we denote uniformly
recurrent as u.r.) word $W$. For any subword $u\sqsubset
W$ lets suppose upper density of occurrence coincides with lower density. The
invariant measure is unique, then.

In what cases $M$ is a torus and $f: M\to M$ is its shift? It means that this
dynamical system has discrete spectrum. Let $W$ be a superword,
obtained by this dynamics. Let $T$ be shift operator. Then
mismatch function between $W$ and $T^n(W)$ satisfies following
conditions:

\begin{enumerate}
    \item There exists sequence $\{n_i\}$ such that
    $\rho(T^{n_i}(W),W)\to 0$.
    \item There exists co-prime arbitrary large $n_i, n_j$ from
    this sequence.
\end{enumerate}

We shall analyse these problems. We start from general
constructions, regarding torus rotation questions, uniqueness of
invariant measure, minimality of dynamical system.

The famous  Sturmian sequences and some of their generalizations
present situation of ``combinatorial paradise''. It provides
patterns for further investigation. Using language of Rauzy graphs
we shall formulate criteria of a superword being generated by interval
exchange transformation. Note that any billiard word with rational
angles can be obtained via such transformations. (Number of
directions of ball is finite and position of ball on the side
together with its direction provides phase point, and Phase space
is union of some intervals.) On the other hand, using language of
Rauzy schemes (obtained from Rauzy graphs by exchanging maximal
sequences of vertices of ingoing and outgoing degree 1 by arches)
we get the criteria of the superword to be morphic.

Concerning shifts of multydimensional torus, there is a beautiful
theory of Rauzy fractals. For more complicated systems one needs
the other patterns for investigation rather then provided by
Sturmian sequences.

The rest of this paper is devoted to dynamic systems connected
with unipotent transformations of torus. This subject was considered
in \cite{ArnoouxMaduit}. Issues related to the study of
sequences, obtained by taking fractional part of values of a
polynomial at integer points lead to the investigation of such
dynamic systems. This problems play an important role in theory of
numbers, theory of information transfer and some other branches
\cite{po,VI,ku_ni}. Note that unipotent transformation of
${\mathbb T}^2$ has the same relation with circle shift as
billiards with arbitrary angles to interval exchange
transformations. Sequences appearing in such billiards
are analysed in \cite{CassHubTroub}.

Inverse problems of symbolic dynamics related to the unipotent
transformation of a torus were studied in paper~\cite{BK}
(Unfortunately it was published only in Russian).

We have to point out the paper \cite{ArnoouxMaduit} obtained
independently from \cite{BK}.
 Let $Q(X)$ be a real polynomial of degree $d\geq 1$ where the
 coefficient of $X^d$ is irrational. Define the difference operator
 $\Delta u_n=u_{n+1}-u_n$ and its iterates
 $\Delta^2=\Delta \circ \Delta,\cdots,\ \Delta^d=\Delta \circ \Delta^{d-1}$.
 The authors establish the theorem according to which the sequence
 $(\Delta^d\lfloor Q(n)\rfloor)_{n\geq 0}$ takes its value on a $2^d$-el\-em\-ent
 alphabet and that
 $$
 p_d(n)=\frac{1}{V(0,1,\cdots, d-1)}\sum_{0\leq k_1<\cdots <k_d\leq n+d-1}V(k_d,\cdots,k_1),
 $$
where
$$V(k_d,\cdots, k_1)=\prod_{1\leq i<j\leq d}(k_j-k_i),$$
 a
Vandermonde determinant. In particular, $p(n)$ depends only on the
degree $d$ of the polynomial $Q$ provided the coefficient of $X^d$
is irrational, $p_2(n)=(n+1)(n+2)(n+3)/6$.

In this situation one have to count the number of parts of torus
division by images of hyperplane. And proof that points in
different regions have different evolutions can be done just by
dimension induction because this system provides more information,
no theory of quasi-invariant sets and factor dynamics required.

The main theorem of this section is

\begin{theorem}        \label{ThMainUnipTorusTransf}
There exists a polynomial $Q(k)$, depending only on the degree of
the polynomial $P$, such that $T(k)=Q(k)$ for all sufficiently large
$k$.
\end{theorem}

\section{Sturmian sequences and their generalizations}

The  Problems (both direct and inverse) related to the rotation of
a circle bring to a class of words which are called {\it
Sturmian words}. Sturmian words are infinite words over a binary
alphabet which contain exactly $n+1$ different subwords (factors)
of length $n$ for any $n\geq 1$.

Sturmian words provide an example of correspondence between
language of dynamical system and combinatorial properties of
superwords.  We shall formulate classical result:

\begin{theorem}[Equivalence theorem (\cite{MorseHedlund},\cite{L1}).]\label{ThEquivClassik}
Let $W$ be an infinite recurrent word over the binary alphabet
$A=\{a,b\}$. The following conditions are equivalent:

\begin{enumerate}
\item The word $W$ is a Sturmian word, i.e., for any $n\geq 1$,
the number of different subwords of length $n$ that occur in $W$
is equal to $T_n(W)=n+1$.

\item The word is not periodic and is {\em balanced}, i.e., any
two subwords $u,v\subset W$ of the same length satisfy the
inequality $||v|_a-|u|_a|\leq 1$, where $|w|_a$ denotes the number
of occurrences of symbol $a$ in the word $w$.

\item The word $W=(w_n)$ is a mechanical word with irrational
$\alpha$, which means that there exist an irrational $\alpha$,
$x_0 \in [0,1]$, and interval $U\subset \mathbb{S}^1$,
$|U|=\alpha$, such that the following condition holds:
$$
w_n=\left\{
\begin{array}{rcl}
   a,&{T_{\alpha}}^n(x_0)\in U\\
   b,&{T_{\alpha}}^n(x_0)\not\in U\\
\end{array}\right.
$$

\item Word  $W$ can be obtained as a limit of the sequence of
finite words $\{w_i\}_{i=1}^\infty$, such that $w_{i+1}$ can be
obtained from $w_{i}$ via substitution of the following type
$a^{k_i}b\to b, a^{k_i+1}b\to a$ or $b^{k_i}a\to a, b^{k_i+1}a\to
b$.

Iff sequence of these substitutions is periodic, then $\alpha$ is
a quadratic irrational.

\end{enumerate}
\end{theorem}

Sturmian words can be considered as theoretical ``paradise'' and
pattern for further investigations. There are several different
ways of generalizing Sturmian words.

First, one can consider {\it balanced} words over an arbitrary
alphabet. Balanced nonperiodic words over an $n$-letter alphabet
were studied in paper~\cite{GR} and later in~\cite{H}. In the
papers~\cite{BC} and \cite{C1}, a dynamical system that generates
an arbitrary nonperiodic balanced word was constructed.

Secondly, generalization may be formulated in terms of the {\it
complexity function}. Complexity function $T_W(n)$ presents the
number of different subwords of length $n$ in the word $W$.
Sturmian words satisfy the relation $T_W(n+1)-T_W(n)=1$ for any
$n\geq 1$. Natural generalizations of Sturmian words are words
with minimal growth, i.e., words over a finite alphabet that
satisfy the relation $T_W(n+1)-T_W(n)=1$ for any $n\geq k$, where
$k$ is a positive integer. Such words were described in terms of
rotation of a circle in paper~\cite{C2}. Note also that words
whose growth function satisfies the relation $\lim_{n \to \infty}
T(n)/n = 1$ were studied in paper~\cite{Ab}.

Words with complexity function $T_W(n)=2n+1$ were studied by
P.~Arnoux and G.~Rauzy (\cite{AR,Ra,Ra1}), words with growth
function $T_W(n)=2n+1$ were analysed by G.~Rote~\cite{R}.
Consideration of the general case of words with linear complexity
function involves the study of words generated by interval
exchange transformations. The problem of description of such words
was posed by Rauzy \cite{Ra1}. Words with linear growth of the
number of subwords were studied by  V.~Berth\'e, S.~Ferenczi, Luca
Q.~Zamboni (\cite{FZ}, \cite{BFZ}). They  investigate
combinatorial sequences related with interval exchange
transformations. See also works of P.Bal\'azi, Z.Mas\'akov\'a,
E.Pelantov\'a (\cite{Ba1},\cite{Ba2}, \cite{BMP}).

\section{Interval exchange transformations}

Sturmian sequences can be obtained via specific rotation of unit
circle. Interval exchange transformation generalizes circle
rotation. G.Rauzy posed the question about description of words
obtained by interval exchange transformation \cite{Ra1}.

S.Ferenci and L.Zamboni \cite{FZ3} obtained the criteria of words
generated by interval exchange transformations with following
condition: trajectory of every break point don't get on an any
break point. In this case obviously a complexity function of the
word is equal to $T(n)=(k-1)n+1$. In fact, this is the answer to
the Rauzy question.

In the papers \cite{BelovChernInt,BelovChernyat}  words generated
by general piecewise-continuous transformation of the interval
were studied. This approach is quite different. The answer to this
question is given in terms of the evolution of the {\it labelled
Rauzy graphs} of the word $W$. The {\it Rauzy graph} of order $k$
(the {\it $k$-graph}) of the word $W$ is the directed graph whose
vertices biuniquely correspond to the factors of length $k$ of the
word $W$ and the vertex $A$ is connected to vertex $B$ by directed
arc iff $W$ has a factor of length $k+1$ such that its first $k$
letters make the subword that corresponds to $A$ and the last $k$
symbols make the subword that corresponds to $B$. By the {\it
follower} of the directed $k$-graph $G$ we call the directed graph
$\Fol(G)$ constructed as follows: the vertices of graph $\Fol(G)$
are in one-to-one correspondence with the arcs of graph $G$ and
there exists an arc from vertex $A$ to vertex $B$ if and only if
the head of the arc $A$ in the graph $G$ is at the notch end of
$B$. The $(k+1)$-graph is a subgraph of the follower of the
$k$-graph; it results from the latter by removing some arcs.
Vertices which are tails of (or heads of) at least two arcs
correspond to {\it special factors}; vertices which are heads and
tails of more than one arc correspond to {\it bispecial factors}.
The sequence of the Rauzy $k$-graphs constitutes the {\it
evolution} of the Rauzy graphs of the word $W$. The Rauzy graph is
said to be {\it labelled} if its arcs are assigned by letters $l$ and
$r$ and some of its vertices (perhaps, none of them) are assigned by
symbol ``-''. The {\it follower} of the labelled Rauzy graph is
the directed graph which is the follower of the latter (considered
a Rauzy graph with the labelling neglected) and whose arcs are
labelled according to the following rule:

\begin{enumerate}
\item Arcs that enter a branching vertex should be labelled by the
same symbols as the arcs that enter any left successor of this
vertex;

\item Arcs that go out of a branching vertex should be labelled by
the same symbols as the arcs that go out of any right successor of
this vertex;

\item If a vertex is labelled by symbol ``--'', then all its right
successors should also be labelled by symbol ``--''.
\end{enumerate}

 The
evolution is said to be {\it correct} if, for all $k\geq 1$, the
following conditions hold when passing from the $k$-graph $G_k$ to
the $(k+1)$-graph $G_{k+1}$ :

\begin{enumerate}
\item Each vertex is an incident to at most two incoming and
outgoing arcs;

\item If the graph contains no vertices corresponding to bispecial
factors, then $G_{n+1}$ coincides with the follower $D(G_n)$;

\item If the vertex that corresponds to a bispecial factor is not
labelled by symbol ``--'', then the arcs that correspond to
forbidden words are chosen among the pairs $lr$ and $rl$;

\item If the vertex is labelled by symbol ``--'', then the arcs to
be deleted should be chosen among the pairs $ll$ or $rr$.
\end{enumerate}

The evolution is said to be {\it asymptotically correct} if this
condition is valid for all $k$ beginning with a certain $k=K$. The
{\it oriented} evolution of the graphs means that there are no
vertices labelled by symbol ``--''.

\begin{theorem}[\cite{BelovChernyat,BelovChernInt}]\label{ThExchBelChern}
 A uniformly recurrent word $W$
 \begin{enumerate}
\item is generated by an interval exchange transformation if and
only if the word is provided with the asymptotically correct
evolution of the labelled Rauzy graphs;

\item is generated by an or\-i\-en\-ta\-ti\-on--pres\-er\-ving
interval exchange transformation if and only if the word is
provided with the asymptotically correct oriented evolution of the
labelled Rauzy graphs.
\end{enumerate}
\end{theorem}

The proof of this theorem consists of two stages. First proves 
that these conditions are sufficient for the word to be generated
by a piecewise-continuous interval transformation. And the second
step is to prove that the sets of uniformly recurrent words
generated by piecewise-continuous interval transformations and by
the interval exchange transformation are equivalent. In order to
do it, there is introduced invariant measure which provides metric:
measure of line segment is its length.

\section{Substitutional sequences}

Important class sequences are so-called {\it substitutional
sequences}. These sequences are invariant under substitution. We
refer to the paper \cite{MauduitSubst}. Fibonacci word is such an
example. Lets consider substitution $\phi: 0\to 001, 1\to 01$. From
symbol $0$ one can get Fibonacci word as superword
$\phi^\infty(0)$. It is Sturmian. Tribonacci word
$\tau^\infty(0)=01020100102010$ can be generated by substitution
$\tau(0)=01, \tau(1)=02,\tau(2)=0$.

 In fact, in \cite{Ra19} G.Rauzy showed that the Tribonacci minimal subshift
 (the shift
orbit closure of t) is a natural coding of a rotation on the
$2$-dimensional torus ${\mathbb T}^2$; i.e., is
measure-theoretically conjugate to an exchange of three fractal
domains on a compact set in ${\mathbb R}^2$. Each domain is
translated by the same vector modulo a lattice.

This is one of the most impressive results. It provides
description of two dimensional spaces. Symbols correspond with the division
of ${\mathbb T}^2$ into fractals, called {\it Rausy fractals}.
Theory of Rauzy fractals was generalized on the so-called {\it
Pisot substitutions}.

This type of dynamical systems provides rich structures with nice
picture in multidimensional case. We have correspondence between
language of arithmetics, dynamical systems and combinatorics.

In other multidimensional systems one can get some information.
Complexity of sequences related with multidimensional systems was
studied in \cite{ArnoouxMaduitShTmrJ}.

\subsection{language of substitutions and Rauzy graphs}

Consider condition (4) in the theorem \ref{ThEquivClassik}.
Sturmian sequence can be obtained as a limit of very specific
substitutions. If it is an invariant under some substitution, then
rotation number is quadratic irrationality. Similar fact is known
for any rotation of circle. However, there exists substitutions
such that eigenvalues of corresponding matrices are algebraic
numbers of an arbitrary degree.

The language of Rauzy graphs provides a bridge between
combinatorial and topological properties in problems regarding
interval exchange transformations. Technique of Rauzy graphs is an
important tool in combinatorics of words.

Rauzy graph $G_k$ of Sturmian sequence has one incoming and one
outgoing branching vertice. When they coincide, (byspecial word
appears) $G_{k+1}$ can be obtained via choosing one of two
possibilities, according to the theorem \ref{ThExchBelChern}. This choice
corresponds with the decomposition of $\alpha$ in chain fraction. If this
choice is made in periodic way $\alpha$ is quadratic irrational.

This fact can be generalized. Let $W$ be a u.r. word. Suppose
behavior of Rauzy graphs is periodic in the same sense, then $W$
is equivalent to a superword invariant under some substitution. In
order to formulate a  theorem one should define what 
``periodicity of events'' in Rauzy graphs means.

{\it Rauzy scheme} of $W$ is a sequence of graphs $\{\Gamma_i\}$
such that every vertice of $\Gamma_i$ is either outgoing or
incoming branching vertice of some order, and corresponds to
subword of $W$. $\Gamma_{i+1}$ can be obtained from $\Gamma_i$ via
exchanging some pathes of length 2 passed through some vertice via
arrows and deleting vertices which are not an endpoints of new
arrows. {\it Rauzy scheme is periodic} if there exists $k>0$ such as for
all sufficiently large $i$ there is isomorphism between $\Gamma_i$
and $\Gamma_{i+k}$ which fits with transmission from $\Gamma_i$ to
$\Gamma_{i+1}$.

\begin{theorem}[\cite{BelovMitrofanov}]
Uniformly recurrent superword is equivalent to image of some
morphism of substitutional invariant sequence iff it has a
periodic Rauzy scheme.
\end{theorem}

Proof of the implication that if $W$ has a periodic Rauzy scheme it is subtitutional is not
so difficult. The main obstacle is the opposite direction.
Substitution can be ``bad'' in the sense that image of letter
$a_i$ can have a ``parasite'' inclusion of image of $a_j$. So this naive
construction fails. We don't know explicit construction.

Let call {\it weight of vertice} of Rauzy scheme as a length of
corresponding word.   In order to prove the theorem, one needs to
establish that if $W$ is uniformly recurrent word stable under
substitution, then ratios of weights of all Rauzy vertices are
bounded. Then one can use J.Cassaigne result \cite{Cassaigne}
saying that if $W$ is uniformly recurrent and $\liminf
T_W(n)/n<\infty$, then $\limsup T_W(n+1)-T_W(n)<\infty$. It
follows from the fact that the number of vertices in Rauzy scheme is
bounded and gets existence of periodic Rauzy scheme from that.
(Condition for morphic uniformly recurrent words $\liminf
T_W(n)/n<\infty$ follows from results of Yu.Pritukin
\cite{MuchnPritSemen}.

Proof of ratios boundness is based on the following fact. If $W$
is a morphic uniformly recurrent word then there exists a constant
$C(W)$ such that for any subword $W$ $u$ occurs
in $v$ for any $v\sqsubset W, |v|>C(W)\cdot|u|$. In order to use
it, one needs to construct sequence of Rauzy schemes in such a way,
that pathes incompatible by inclusion correspond to the  words
with the same property.

This theorem implies Vershik-Lifshic theorem \cite{Vershik2,VL} of
periodicity of Bratelli diagrams of Markov compact corresponding
to the substitutinal systems.

A {\it Bratteli diagram} $(V,E)$ is a countable collection $V$ of
finite vertex sets, $V =\{V_n\}_{n=1}^\infty$ and a countable
collection $E$ of finite edge sets $E =\{E_n\}_{n=1}^\infty$,
along with functions $s: E_n\to V_{n-1}$ and $r: E_n\to V_n$ such
that (i) $V_0 = \{\nu_0\}$, (ii) $s: E_n\to V_{n-1}$ and $r:
E_n\to V_n$ are onto for all n. We view $(V, E)$ as a directed
graph where an edge $e\in E_n$ connects the source vertex $s(e)\in
V_n$ to the range vertex $r(e)\in V_{n+l}$. {\it Periodicity} of
Bratelli diagram means that for some $k$ there exists a pair of
mappings from $V_n$ to $V_{n+k}$ and from $E_n$ to $E_{n+k}$
preserving functions $s$ and $r$. With Bratelli diagram one can
associate topological dynamics. Details can be founded in
\cite{HPS}.

The proof of the Vershik-Lifshic theorem uses straightforward
construction. Consider image of
$\varphi^{(n)}(a)=(\varphi^{(n-2)})(\varphi^{(2)}(a))$ for some
letter $a$. It consists of blocks corresponding application of
$\varphi^{(n-2)}$ to the letters of $\varphi^{(2)}(a)$ and also
can be decomposed to the blocks corresponding application of
$\varphi^{(n-1)}$ to the letters of $\varphi(a)$. Finite sets
forming Bratelli diagrams, consists of sequences of pairs (block,
its position in bigger block), they corresponds some subwords of
$W$. On these sets relations of being left and right neighbors
inside bigger block can be naturally posed. The first and last
occurrence in the bigger block needs special attention, because
bigger block itself may have position and can be preceded or
followed by another bigger block. Details can be found in
\cite{Vershik2,VL}. See also \cite{LifshicDisser}.

\section{Main constructions and definitions}

\subsection{Complexity function, special factors, and uniformly recurrent
words} \label{SbScCmplUnirec}

In this section we define the basic notions of combinatorics of
words. By $L$ we denote a finite alphabet, i.e., a nonempty set of
elements (symbols). We use the notation $A^+$ for the set of all
finite sequences of symbols or {\it words}.

A finite word can always be uniquely represented in the form $w =
w_1 \cdots w_n$, where $w_i \in A$, $1\leq i \leq n$. The number
$n$ is called the {\it length} of word $w$; it is denoted by
$|w|$.

The set $A^+$ of all finite words over $A$ is a simple semigroup
with concatenation as semigroup operation.

If element $\Lambda$ (the empty word) is included in the set of
words, then this is actually the free monoid $A^*$ over $A$. By
definition the length of the empty word is $|\Lambda|=0$.

A word $u$ is a {\it subword} (or {\it factor}) of a word $w$ if
there exist words $p,q\in A^+$ such that $w = puq$.

Denote the set of all factors (both finite and infinite) of a word
$W$ by $F(W)$. Two infinite words $W$ and $V$ over alphabet $A$
are said to be {\it equivalent} if $F(W)=F(V)$.

The {\it beginning $w^b$} of the word $w$ is a sequence
$$x,f(x),\ldots,f^{(n)}(x),{\-} n~\in~{\nit}.$$

We say that symbol $a\in A$ is a {\it left} (accordingly, {\it
right}) {\it extension} of factor $v$ if $av$ (accordingly, $va$)
belongs to $F(W)$. A subword $v$ is called a {\it left}
(accordingly, {\it right}) {\it special factor} if it possesses at
least two left (right) extensions. A subword $v$ is said to be
{\it bispecial} if it is both a left and right special factor at
the same time. The number of different left (right) extensions of
a subword is called the {\it left (right) valence} of this
subword.

A word $W$ is said to be {\it recurrent} if each its factor occurs
in infinitely many times (in the case of a doubly-infinite
word, each factor occurs infinitely many times in both
directions). A word $W$ is said to be {\it uniformly recurrent} or
({\it u.r word}) if it is recurrent and, for each of its factor
$v$, there exists a positive integer $N(v)$ such that, for any
subword $u$ of length at least $N(v)$ of the word $W$, factor $v$
occurs in $u$ as a subword.

Below we formulate several theorems about u.r words, which will be
needed later. The proof of these theorems can be found in
monograph~\cite{BBL}.

\begin{theorem}
The following two properties of an infinite word $W$ are
equivalent:

a) For any $k$ there exists $N(k)$ such that any segment of length
$k$ of the word $W$ occurs in any segment of length $N(k)$ of the
word $W$;

b) If all finite factors of a word $V$ are at the same time finite
factors of a word $W$, then all finite factors of the word $W$ are
also finite factors of the word $V$.
\end{theorem}

\begin{theorem}\label{ThUnifRec}
Let $W$ be an infinite word. Then there exists a uniformly
recurrent word $\widehat{W}$ all of whose factors are factors of
$W$.
\end{theorem}

One can consider the action of the shift operator $\tau$ on the
set of infinite words. {\it The Hamming distance} between words
$W_{1}$ and $W_{2}$ is the quantity $d(W_{1},W_{2}) =\sum_{n\in
{\mathbb Z}}\lambda _{n}2^{-|n|}$, where $\lambda _{n}= 0$ if
symbols at the $n$-th positions of the words are the same and
$\lambda _{n}= 1$, otherwise.

An {\it invariant subset} is a subset of the set of all infinite
words which is invariant under the action of $\tau$. A {\it
minimal closed invariant set}, or briefly, m.c.i.s, is a closed
(with respect to the Hamming metric introduced above) invariant
subset which is nonempty and contains no closed invariant subsets
except for itself and the empty subset.

\begin{theorem}[Properties of closed invariant sets]\label{RecEq}
The following properties of a superword $W$ are equivalent:
\begin{enumerate}
\item $W$ is a uniformly recurrent word;

\item The closed orbit of $W$ is minimal and is a m.c.i.s.
\end{enumerate}
\end{theorem}

\begin{theorem}\label{RecCont}
Let $W$ be a uniformly recurrent nonperiodic infinite word. Then
\begin{enumerate}
\item All the words that are equivalent to $W$ are u.r. words; the
set of such words in uncountable;

\item There exist distinct u.r. words $W_1 \neq W_2$ which are
equivalent to the given word and can be written as $W_1=U V_1$,
$W_2=U V_2$, where $U$ is a left-infinite word and $V_1\neq V_2$
are right-infinite words.
\end{enumerate}
\end{theorem}


\section{Unipotent dynamics on a torus}

\subsection{Essential evolution of points.}

Let $f:M~\to~M$ be a continuous map on the space $M$ and $U
\subset M$ be a subset. The starting point $x$ determines a binary
word $w$ describing the evolution as above: $w_n=1$, if
$f^{(n)}(x)\in U$ and $w_n=0$ if $f^{(n)}(x)\not\in U$.  We assume
that $U$ is an open set, $\mes(\partial U) = 0$ and $M$ is a
compact metric space.


\begin{definition}

A finite word $v^f$ is said to be an {\em essential finite
evolution of a point $x^{\ast}$}, if every neighborhood of a point
$x^{\ast}$ contains an open set $V$, so that for all $x\in V$
holds $v_x^b=v^f$.

An infinite word $w$ is said to be an {\em essential (infinite)
evolution of a point $x^{\ast}$}, if every initial subword is an
essential finite evolution of a point $x^{\ast}$.

\end{definition}

The word ``evolution'' will further denote ``essential
evolution''.

\begin{proposition}
Let $v$ be a finite word, then the set of points with fixed finite
evolution is closed (i.e. consists all its limit points)

The similar statement holds for an infinite word $w$. (The
intersection of any family of closed sets is closed)
\end{proposition}

We shall not use the next proposition although it has its own
interest:

\begin{proposition}
Let $(M,f,U,x)$ be a dynamical system without closed invariant
subsets, and different points have different evolutions. Let
$(\hat{M},s,U,x)$ be corresponding symbolical dynamical system,
i.e. set of all superwords with Tikhonov topology. Then $M$ is
naturally isomorphic to factor $\hat{M}$ by spaces consisting of sets
of superwords which are essential evolution of one point from $M$,
isomorphism induces natural isomorphisms of dynamical systems.
\end{proposition}

\subsection{Morphism of dynamics}
\begin{definition}

A morphism of two dynamics $G:(M_1,f_1) \to (M_2,f_2)$ is a
continuous map, such that the diagram
\begin{center}
$$
\begin{array}{clcll}
{}&M_1& \longrightarrow^{\hskip -11pt g}& M_2&{}\cr {}& f_1
\downarrow&{}&\downarrow f_2 \cr {}& M_1& \longrightarrow^{\hskip
-11pt g}& M_2
\end{array}
$$
\end{center}
is commutative.

\end{definition}

The notions of {\it epimorphism, monomorphism} and {\it
isomorphism} are defined in the natural way.

The {\it factor-dynamics} is naturally defined on the quotient
topology iff $f$ permutes the equivalence classes of the map $f$.
Note that the inverse images of points under morphisms are closed.

\begin{definition}
A set $V$ is {\em irreducible}, if its closure is not an inverse
image of a closed set under any morphism (except monomorphism) and
reducible otherwise.
\end{definition}

\begin{theorem}\label{th1}
Let a set $U$ be  irreducible, then
\begin{enumerate}
\item Different points have different evolutions.

\item For any $\varepsilon > 0$ there exists $N(\varepsilon)$,
such that every two words of length $N(\varepsilon)$,
corresponding to the initial points on distance greater than
$\varepsilon$ are different.
\end{enumerate}
\end{theorem}

\begin{proof}
\ \
\begin{enumerate}
\item Classes of points with the same evolution are closed and the
map $f$ permutes them.

\item The proof of this point follows from the proof for the
previous one by contradiction and transition to limit.
\end{enumerate}
\end{proof}

\subsection{Quasi-invariant sets}

\begin{definition}
A dynamics is said to be {\em minimal}
 if $M$ does not contain closed invariant sets
apart from  $M$ and $\O$ and {\em irreducible}, if it does not
contain proper Quasi-invariant sets.
\end{definition}

\begin{definition}
A closed set $N$ is {\em quasi-invariant} if for every two points
$A$ and $B$ and a convergent sequence $f^{(n_i)}(A) \to C$ such
that $C \in N$  for $n_i \to \infty$ every limit point of the
sequence $f^{(n_i)}(B) \in N$.
\end{definition}

\begin{proposition}\label{st2}
\ \
\begin{enumerate}
\item Every closed invariant set is quasi-invariant.

\item   A partition of
qu\-a\-si-in\-va\-ri\-ant sets corresponds to each factor-dynamics and conversely.

\item The image of a quasi-invariant set is quasi-invariant.

\item
If $M$ is closed invariant set, than one qu\-a\-si-in\-va\-ri\-ant set
uniquely determines the fac\-tor--dy\-na\-mic.

\item The set of points with fixed evolution is quasi-invariant.
\end{enumerate}
\end{proposition}

\begin{definition}
{\em $A-B$ cloud (or 0-cloud) with center $A$, generated by point
$B$} is the set of conditional limit points $f^{(n_i)}(B)$ under
the condition $f^{(n_i)}(A) \to A$.

{\em $A-B$ $k$-cloud with center $A$, generated by point $B$} is
the closure of the union of set of conditional limit points
$f^{(n_i)}(B)$ under the condition $f^{(n_i)}(A) \to A^{\ast}$
where $A^{\ast} \in A-B$ $(k-1)$.

{\em $A-B$ $k$-cloud with center $A$, generated by point $B$} is
the closure of the union of the set of $A-B$ $k$-clouds, $k \in
N$.
\end{definition}

Note that  $A-B$ $k$-cloud is closed.

\begin{proposition}\label{PRCloud}

a) The image of $A-B$ $k$-cloud under the $l$-th iteration is
$f^{(l)}(A)-f^{(l)}(B)$-clouds.

b) If $A_n \to A$, $B_n \to B$ then $\rho(A_n-B_n,A-B) \to 0$

\end{proposition}

Denote $A-B$-cloud by $l_0$, by  $L_{i+1}$ lets denote the closure of the
union of all $A_i-B_i$-clouds, for which $A_i,\ B_i \in L_i$. Assume
$L_{AB}= \bigcup L_i$. Factorization, generated by $L_{AB}$ is the
weakest factorization, that glue points $A$ and $B$.

\subsection{Unipotent dynamics on a torus}

Problems, connected with the study of the behavior of fractional
parts of values of a polynomial at integer points, are in fact the
classical problems of symbolic dynamics. Let $P(n)$ be a polynomial
of degree $m+1$ with irrational coefficient $a_{m+1}$ of the
highest degree. Define a sequence of polynomials $P_k(n)\ \
k=0,\dots,m$ in the following way:

\begin{eqnarray}\label{eq1}
\left\{
\begin{array}{ccc}
P_m(n)=P(n), \nonumber \\
P_{m-1}(n)=P_m(n+1)-P_m(n), \nonumber \\
\ldots \\
P_{i-1}(n)=P_i(n+1)-P_i(n), \nonumber \\
\ldots \nonumber
\end{array}
\right.
\end{eqnarray}

From these formulas it follows, that $P_0(n)=n!a_{m+1}$ is
irrational. Put $\varepsilon=P_0(n)$, $x_i(n)=\{P_i(n)\}$ and
$x_i^{\prime}(n)=x_i(n+1)$, then from (\ref{eq1}) we obtain the
following dynamical system:

\begin{eqnarray}\label{eq2}
\left\{
\begin{array}{lccc}
x_m^{\prime}\ =\ x_{m}+x_{m-1}\ \ \mod\ 1\\
x_{m-1}^{\prime}\ =\ x_{m-1}+x_{m-2}\ \ \mod\ 1\\
\cdots \\
x_1^{\prime}\ =\ x_{1}+\varepsilon\ \ \mod\ 1,
\end{array}
\right.
\end{eqnarray}

where $\varepsilon$ is irrational as $\varepsilon=n!a_{m+1}$. The
condition $[2\{P(n)\}]=0$ turns to the condition $0 \le x_m(n) <
1/2$.

Consequently the vector $(x_1^{\prime},\ldots,x_m^{\prime})$ is
obtained from $(x_1,\cdots,x_m)$ by unipotent transformation
(transformation, corresponding to a linear transformation with
unitary eigenvalues).

Images of hyperplanes  $x_m=0$ and $x_m=1/2$ divide the space on
polyhedra. The same word of length $k$ corresponds with the points of the same polyhedron.


\subsection{Mismatch function}

\begin{definition}
A function $\rho$ is said to be a mismatch function of words $w$
and $v$ and is defined in the following way:
\begin{eqnarray*}
\rho(i)=\left\{
\begin{array}{ccc}
0,\ \ if \ \ w_i=v_i,
\\ 1,\ \ if \ \ w_i\not=v_i,
\end{array}
\right.
\end{eqnarray*}

A density of mismatch of $\rho(w,v)$ of words $w$ and $v$ is defined
by the formula
\begin{displaymath}
\rho(w,v)=\lim_{i\to\infty} \frac{\sum_{j=1}^{i} \rho(j)}{i}
\end{displaymath}
\end{definition}

\begin{theorem}
Let $w$ � $v$ be two different evolutions of the point $x_0\in T$.
Then $\rho(w,v)=0$.
\end{theorem}
\begin{proof}

From lemma Wail  \cite{Wail} it follows, that the orbit of any
point $x_0\in T$ is everywhere dense and evenly distributed. From the
continuity of map $f$, the condition $\mes(\partial U) = 0$ and
the definition of evolution it follows,  that  $w_n=v_n$ if
$f^{(n)}(x_0)\not\in
\partial U$. The proof of the theorem follows from these statements.
\end{proof}

\begin{theorem}
Let points $x$ and $x^{\ast}$ be different and have different
evolutions. Then the density of mismatch $\rho(w_x,w_{x^{\ast}})$
is defined and is greater than $0$.
\end{theorem}

\begin{proof}
We can assume that words  $w_x$ and $w_{x^{\ast}}$ differ in the
first position. Consider the direct product  $T\times T$. We
divide the set of points into two classes: with the same current
position and different current position. Let $O \subset T\times T$
be the set of different pairs, then the orbit of pair
$(x,x^{\ast})$ lies in $O$. Closed orbit of any pair is a torus, a
minimal closed invariant set, on which the dynamics on torus is
realized.  Let $\rho$ be the volume of the intersection, then
$\rho \not= 0$ and $\rho$ is the density.

\end{proof}

\subsection{Description of torus factor-dynamics.}

We'll consider the dynamics that don't glue the coordinate  $x_n$.
Consider the $k$-th iteration of the transformation of torus:
\begin{eqnarray}\label{eq3}
\left \{
\begin{array}{ccc}
x_m^{(k)}=x_{m}+C_k^1x_{m-1}+\ldots+C_k^{m}\varepsilon\ \ \mod 1\\
\cdots \\
x_i^{(k)}=x_{i}+C_k^1x_{i-1}+\ldots+C_k^{i}\varepsilon\ \ \mod 1\\
\cdots \\
x_1^{(k)}=x_1+C_k^1\varepsilon\ \ \mod 1\\
\end{array}
\right.
\end{eqnarray}

If the points $A$ $(x_1,\dots,x_m)$ and $B$ $(x_1+\Delta
x_1,\dots,x_m+\Delta x_m)$ belong to the set $M^{\prime}$, then
$f^{(k)}(A)$ and $f^{(k)}(B)$ also belong to the same set. From
(\ref{eq3}) it follows:

\begin{eqnarray}\label{eq4}
\left \{
\begin{array}{ccccccc}
\Delta x_m^{(k)}=\Delta x_{m}+C_k^1\Delta
x_{m-1}+\ldots+C_{k}^{m-1}\Delta
x_1\ \ \mod\ 1\\
\cdots \\
\Delta x_i^{(k)}=\Delta x_{i}+C_k^1\Delta
x_{i-1}+\ldots+C_{k}^{i-1}\Delta x_1\ \ \mod\ 1\\
\cdots \\
\Delta x_1^{(k)}=\Delta x_1\ \ \mod\ 1\\
\end{array}
\right.
\end{eqnarray}
\begin{proposition}
Closed orbit of any pair of points is a torus $T^{\prime}$, a
minimal closed invariant set, on which the dynamics of torus is
realized.
\end{proposition}


%
%

\begin{remark} For two-dimensional case there exists an area $U$
with analytic boundary and a starting point $x^{\ast}$ that have
infinitely many different evolutions, that differ in an infinite
number of positions.
\end{remark}

%

\begin{proposition}\label{st4}
Let $1$, $\varepsilon$, $\Delta x_i$ be linearly independent over
$\mathbb{Q}$. Then $A-B$ cloud contains all points which first $i-1$
coincide with coordinates of point $B$.
\end{proposition}

\begin{proof}
We may assume that $\Delta x_j$ is rational when $j<i$ and let
$\Delta x_j=p_j/q_j$ be the presentation $\Delta x_j$ as an
irreducible fraction and $k$ is divisible by the product
$m!\prod_{j=1}^i q_j$, then $x_j^{(kl)}=x_j$ when $j<i$ and systems
(\ref{eq3}) and (\ref{eq4}) can be rewritten as
\begin{eqnarray}\label{eq5}
\left \{
\begin{array}{ccc}
x_m^{(kl)}=x_{m}+C_{kl}^1x_{m-1}+\ldots+C_{kl}^{m}\varepsilon\ \
\mod\ 1\\
\cdots \\
x_j^{(kl)}=x_{j}+C_{kl}^1x_{j-1}+\ldots+C_{kl}^{j}\varepsilon\ \
\mod\ 1\\
\cdots \\
x_1^{(kl)}=x_1+C_{kl}^1\varepsilon\ \ \mod\ 1\\
\Delta x_m^{(kl)}=\Delta x_{m}+C_{kl}^1\Delta
x_{m-1}+\ldots+C_{kl}^{m-1}\Delta x_1\ \ \mod\ 1\\
\cdots \\
\Delta x_j^{(kl)}=\Delta x_{j}+C_{kl}^1\Delta
x_{j-1}+\ldots+C_{kl}^{j-1}\Delta x_1\ \ \mod\ 1 (j \ge i) \\
\Delta x_j^{(kl)}=\Delta x_j\ \ \mod\ 1 (j <i) \\
\end{array}
\right.
\end{eqnarray}

Vector $(x_1^{(kl)}, \ \ldots , \ x_m^{(kl)},\Delta x_i^{(kl)}, \
\ldots , \ \Delta x_j^{(kl)})$ by lemma Wail \cite{Wail} is
everywhere dense in the torus of dimension $2m-i+1$, then from
definition of
 $A-B$ cloud it follows that it contains all points, which first
$i-1$ coordinates coincide with coordinates of the point $B$, and
other points can be chosen arbitrarily.
\end{proof}

\begin{proposition}\label{st5}
Let $\Delta x_i$ be an irrational number. Then there exists a point
$B_n$, the evolution of which coincide with evolution of points  $A$
and $B$, and for which it holds:
$$
\Delta x_i^{B_n} = n \Delta x_i^{B}.
$$
\end{proposition}

\begin{proof}
Lets choose points $A$ and $B$ as $B_0$ and $B_1$ respectively. By the
method of mathematical induction we'll assume that point $B_k$ is
already built, and as for $B_{k+1}$ it suffices to take any
conditionally limit point of a sequence $f^{(n_i)}(B_k)$, under the
condition $f^{(n_i)}(B_{k-1}) \to B_k$. Note, that the point $B_k
\in A-B$ $k$-cloud.

\end{proof}

\begin{proposition} \label{PrAlmFinal}
Let $\Delta x_i$ be an irrational number. Then there exists a point
$B_{\delta}$, the evolution of which coincide with evolutions of
points  $A$ and $B$, for which it holds:
$$
\Delta x_i^{B_{\delta}} =  \delta, 0 \le \delta \le 1.
$$
\end{proposition}

This fact follows directly from the proposition \ref{st5} and the
fact that the set of point with fixed evolution is closed.

Thus the case when $\Delta x_i$ is irrational reduces to the case
when $1$, $\varepsilon$, $\Delta x_i$ are linearly independent over
$\mathbb{Q}$.

The case when all  $\Delta x_j$ are rational leads to
factor-dynamics where $i$ edge of torus divides in
$M_i=\prod_{j=1}^{i}m_j$ parts, where $m_j$ are arbitrary natural
numbers and points of the torus $x=x^{\ast}$ are identified when
for all $1 \le j \le m$ holds: $M_jx_j=M_jx^{\ast}_j$.

Description of Quasi-invariant sets follows from the proposition
\ref{PrAlmFinal}.

\begin{theorem}[Description of Quasi-invariant sets]\label{Thco1}
Quasi-invariant set is a shift of abelian group which transforms
into themselves under translations by $1/M_i$ along $i$-coordinate
or under all translations along coordinates with the number,
greater than some fixed number.
\end{theorem}

The theorem provides description of all possible factor-dynamics.

\begin{corollary}\label{co1}
The class of closed reducible sets consists of sets, which transform
into themselves under translations by $1/M_i$ along $i$-coordinate
or under all translations along coordinates with the number, greater
than some fixed number. The set $0 \le x_m < 1/2$ is irreducible.
\end{corollary}

\subsection{Proof of the theorem \ref{ThMainUnipTorusTransf}}

This theorem follows from the next proposition:

\begin{proposition}\label{tt}
Consider the dynamics of the torus, given by equation (\ref{eq3})
and let irreducible set $U$ is given by the condition:

\begin{displaymath}
0 \le x_m \le 1/2,
\end{displaymath}

then there exists such $L(\varepsilon)$ that points with the same
finite evolution of length $L(\varepsilon)$ divide torus on closed
convex polyhedra, and different polyhedra correspond with the different
evolutions.
\end{proposition}

\begin{proof}
The set $U$ is irreducible, so by theorem \ref{th1} there exists
$\delta^{\prime}$ such that $N(\delta^{\prime})$ evolutions of points
at a distance greater than  $\delta^{\prime}$ are different. It is
obvious, that $N$-evolutions of all interior points of the
polyhedron, obtained after $i$-th iteration are the same. We call it
the evolution of the polyhedron.

Consider parts of the partition, corresponding to words of length
$N(\delta^{\prime})$, then each polyhedron of the partition can't
intersect more than one hyperface from $n+1$ family of planes
$f^{n}[x_1=q/2,\ q \in N$, otherwise we can choose two points, that
belong to this polyhedron and are divided by two planes from $n+1$
family.

 Assuming $n+1$ family to be the first and turning the time back we get the
 situation when $N(\delta^{\prime})$ evolution of these points
 doesn't allow to distinguish them, and this is impossible. Polyhedra
 formed by $n+1$ family of planes with the same evolution can't have
 common points. Let $\delta^{\ast}$ be the minimal distance between
 polyhedra with the same evolutions, then by  putting
 $L(\varepsilon) = \max (N(\delta^{\prime}), \ N(\delta^{\ast}))$
 we obtain the number of iterations, from which certainly achieve
 the equality between number of words of length $k$ and the number
 of polyhedra in which the torus divides $k$ families of planes
 $f^{i}[x_1=q/2,$ $q=0,\ldots,k-1$. When $(\varepsilon)$ is
 irrational the intersection of more than $m$ hyperplanes is empty
 and there exists a one-to-one correspondence between points of
 intersection of $m$ non-parallel planes and polyhedra of partition.
 By computing the number of points in the intersection of $k$
 families of planes $f^{i}[x_1=q/2,$ $q=0,\ldots,k-1$ we obtain a
 polynomial, which defines the number of parts and consequently the
 number of subwords from some moment.
\end{proof}

 The number of points of the intersection of hyperplanes $Q(k)$
and consequently the number of subwords $T(k)$ $(k \ge K)$ of
length $k$ can be calculated by the formula:

\begin{displaymath}
Q(k)= \sum_{0\le k_1 < \ldots < k_m \le k} \left|
\begin{array}{cccc}
1 & {k_m \choose 1} & \ldots & {k_m \choose m} \\
\ & \ & \ldots & \ \\
1 & {k_1 \choose 1} & \ldots & {k_1 \choose m}
\end{array} \right|,\quad \deg Q(k) = \frac{m(m+1)}{2}.
\end{displaymath}

For every $K$ there exists $P$ and $k_0 \geqslant K$ such that the
equality $T(k)=Q(k)$ holds for $k \geqslant k_0$ and doesn't hold
for $k_0$.




\end{document}